\documentclass[12pt,twoside]{article}

\usepackage{a4,graphicx,amsmath,amsfonts,amssymb,mathrsfs,bbm}

\newcommand{\real}{\mathbbm{R}}
\newcommand{\complex}{\mathbbm{C}}

\newcommand{\ltwo}{\mathscr{L}_2[0,\infty)}
\newcommand{\htwo}{\mathscr{H}_2}
\renewcommand{\top}{{\rm T}}

\newtheorem{definition}{Definition}
\newtheorem{theorem}{Theorem}

\newtheorem{lemma}{Lemma}

\begin{document}

\thispagestyle{plain}

\begin{center}
  {\bf \Large Stability preservation in \\[0.5ex]
    Galerkin-type projection-based \\[1.5ex]
    model order reduction}


\vspace{10mm}

{\large Roland~Pulch} \\[1ex]
{\small Institute of Mathematics and Computer Science, \\
Ernst-Moritz-Arndt-Universit\"at Greifswald, \\
Walther-Rathenau-Str.~47, 17489 Greifswald, Germany. \\
Email: {\tt roland.pulch@uni-greifswald.de}}

\end{center}

\bigskip\bigskip


\begin{center}
{Abstract}

\begin{tabular}{p{13cm}}
We consider linear dynamical systems consisting of
ordinary differential equations with high dimensionality.
The aim of model order reduction is to construct an
approximating system of a much lower dimension.
Therein, the reduced system may be unstable, even though
the original system is asymptotically stable.
We focus on projection-based model order reduction of Galerkin-type.
A transformation of the original system guarantees an asymptotically stable
reduced system.
This transformation requires the numerical solution of a high-dimensional
Lyapunov equation.
We specify an approximation of the solution,
which allows for an efficient iterative treatment of
the Lyapunov equation under a certain assumption.
Furthermore, we generalise this strategy to preserve the asymptotic stability
of stationary solutions in model order reduction of
nonlinear dynamical systems.
Numerical results for high-dimensional examples confirm the computational
feasibility of the stability-preserving approach.

\bigskip


\bigskip

\end{tabular}
\end{center}

\clearpage

\markboth{\em R.~Pulch}{\em Stability in Galerkin-type MOR}


\section{Introduction}
The mathematical modelling of problems from science and engineering
often yields dynamical systems.
The increasing complexity of industrial applications causes high-dimensional
systems by electronic design automation, for example.
Thus a numerical simulation of the model may become too costly.
Methods of model order reduction (MOR) are required to decrease
the dimensionality of the dynamical systems,
see~\cite{antoulas,benner-mehrmann,schilders}.

We consider linear implicit systems of
ordinary differential equations (ODEs),
which are asymptotically stable. 
Projection-based MOR yields linear ODEs of a lower dimensionality.
However, the reduced system may be unstable and thus useless.
Firstly, some solutions become unbounded in the time domain. 
Secondly, error bounds, which follow from the transfer functions
in the frequency domain, are not valid any more.
Hence stability-preserving MOR methods are necessary to guarantee
adequate reduced systems.

The method of balanced truncation, see~\cite{gugercin-antoulas},
always yields stable reduced systems, while the computational effort
is often relatively large.
Krylov subspace techniques, see~\cite{freund}, are cheaper,
whereas stability can easily be lost.
A stability-preservation of a Krylov subspace approach is given by
special assumptions and methods in~\cite{ionescu}.
The stability property can also be satisfied by a post-processing
using the poles of the transfer function, see~\cite{bai-freund}.

We examine projection-based MOR of Galerkin-type, where each scheme
is defined by a single orthogonal projection matrix.
Prominent methods are the one-sided Arnoldi algorithm and
the proper orthogonal decomposition (POD), for example.
Prajna~\cite{prajna} introduced an approach to guarantee the stability in
such an MOR of a nonlinear system of ODEs by a basis transformation
in both the state space and the image space.
This technique was also applied to a stochastic Galerkin projection
in~\cite{pulch-augustin}.
Casta{\~n}{\'e}~Selga et al.~\cite{castane-selga} investigated a
stabilisation of linear systems of ODEs by a transformation
in the image space only.
The main effort consists in solving a single Lyapunov equation,
whose efficient numerical solution is critical.
The high dimensionality excludes direct methods and thus
approximate methods have to be used.
The alternating direction implicit (ADI) algorithm,
see~\cite{li-white,penzl}, requires an input matrix
in the form of a symmetric low-rank factorisation.

We perform a specific ansatz for an approximate solution of the
Lyapunov equation.
Our strategy yields an alternative Lyapunov equation, where a
symmetric factorisation can be chosen as input matrix.
Its rank depends on the signs of eigenvalues in a symmetric part.
Consequently, the ADI algorithm is applicable with a small number of
iteration steps in the case of low ranks.
Moreover, our method guarantees a non-singular and often well-conditioned
mass matrix in the reduced system.
The technique can also be employed to preserve the asymptotic stability
of stationary solutions in an MOR of nonlinear  ODEs.

The paper is organised as follows.
We describe the class of MOR methods as well as the stability-preserving
transformation in Section~\ref{sec:problem-def}.
Our numerical method is derived and analysed in Section~\ref{sec:method}.
We demonstrate the applicability of this approach to nonlinear
dynamical systems in Section~\ref{sec:nonlinear}.
Finally, Section~\ref{sec:example} presents numerical results of
illustrative examples.


\section{Projection-based reduction and stability}
\label{sec:problem-def}
The projection-based MOR and the stability problem are defined
in this section.

\subsection{Linear dynamical systems and stability}
Let a linear dynamical system be given in the form of implicit ODEs
\begin{equation} \label{linear-ode}
  \begin{array}{rcl}
    E \dot{x}(t) & = & A x(t) + B u(t) , \\
    y(t) & = & C x(t) \\
  \end{array}
\end{equation}
with state variables $x : [0,\infty) \rightarrow \real^n$,
inputs $u : [0,\infty) \rightarrow \real^{n_{\rm in}}$ and
outputs $y : [0,\infty) \rightarrow \real^{n_{\rm out}}$.
The system includes constant matrices $A,E \in \real^{n \times n}$,
$B \in \real^{n \times n_{\rm in}}$ and $C \in \real^{n_{\rm out} \times n}$.
We assume that the mass matrix~$E$ is non-singular.
Initial values $x(0)=x_0$ are predetermined.

In the frequency domain, a transfer function completely describes
the input-output behaviour of the system~(\ref{linear-ode}),
see~\cite{antoulas}.
This transfer function
$H : \complex \backslash \Sigma \rightarrow \complex^{n_{\rm out} \times n_{\rm in}}$
reads as
\begin{equation} \label{transfer}
  H(s) = C ( s E - A )^{-1} B
  \qquad \mbox{for} \;\; s \in \complex \backslash \Sigma .
\end{equation}
The mapping~(\ref{transfer}) is a rational function with a
finite set of poles~$\Sigma \subset \complex$.
The magnitude of a transfer function can be characterised by norms
in Hardy spaces.
The $\htwo$-norm is defined by, see~\cite[p.~92]{shmaliy},
\begin{equation} \label{h2-norm}
  \left\| H \right\|_{\htwo} =
  \sqrt{ \frac{1}{2\pi} \int_{-\infty}^{+\infty}
    \left\| H({\rm i}\omega) \right\|_{\rm F}^2 \; {\rm d}\omega } 
\end{equation}
with ${\rm i} = \sqrt{-1}$ and the Frobenius (matrix) norm
$\| \cdot \|_{\rm F}$ provided that the integral exists.

The stability issues of the system~(\ref{linear-ode}) are independent of
the definition of inputs and outputs. 
To discuss the stability, we recall some general properties of matrices.

\begin{definition} \label{def:spectral}
  Let $A \in \real^{n \times n}$ and $\lambda_1,\ldots,\lambda_n \in \complex$
  be its eigenvalues.
  The {\em spectral abscissa} of the matrix~$A$ is the real number
  $$ \alpha (A) = \max
  \left\{ {\rm Re}(\lambda_1) , \ldots , {\rm Re} (\lambda_n) \right\} . $$
\end{definition}

\begin{definition} \label{def:stable-matrix}
  A matrix $A \in \real^{n \times n}$ is called a {\em stable matrix},
  if its spectral abscissa satisfies $\alpha(A) < 0$.
\end{definition}

\begin{definition} \label{def:diss-matrix}
  A matrix $A \in \real^{n \times n}$ is called a {\em dissipative matrix},
  if its {\em symmetric part}
  $A_{\rm sym} = A + A^\top$ is negative definite.
\end{definition}

The definition of dissipativity is in agreement to~\cite[p.~62]{gil}.
A dissipative matrix is also stable.
Vice versa, a stable matrix is not dissipative in general.
The linear dynamical system~(\ref{linear-ode}) is asymptotically stable,
if and only if each eigenvalue $\lambda \in \complex$ satisfying
\begin{equation} \label{gen-ev}
  \det \left( \lambda E - A \right) = 0
\end{equation}
has a strictly negative real part,
see~\cite[p.~376]{braun}.
Equivalently, $E^{-1}A$ is a stable matrix.
The asymptotic stability guarantees the existence of the
integral in~(\ref{h2-norm}).
If the eigenvalues of the problem~(\ref{gen-ev}) all have a non-positive
real part and a real part zero appears, then Lyapunov stability may still
be satisfied.
However, we consider this case also as a loss of stability,
because the asymptotic stability is not valid any more.
Likewise, the linear dynamical system~(\ref{linear-ode}) is called
dissipative, if and only if $E^{-1}A$ is a dissipative matrix.


\subsection{Projection-based model order reduction}
\label{sec:mor}
We assume that the linear dynamical system~(\ref{linear-ode}) exhibits
a huge dimensionality~$n$.
Thus the involved matrices~$A$ and~$E$ typically are sparse.
The purpose of MOR is to decrease the complexity.
An alternative linear dynamical system
\begin{equation} \label{system-reduced}
  \begin{array}{rcl}
    \bar{E} \dot{\bar{x}}(t) & = & \bar{A} \bar{x}(t) + \bar{B} u(t) , \\
    \bar{y}(t) & = & \bar{C} \bar{x}(t) \\
  \end{array}
\end{equation}
has to be constructed with state variables
$\bar{x} : [0,\infty) \rightarrow \real^r$
and the matrices $\bar{A},\bar{E} \in \real^{r \times r}$,
$\bar{B} \in \real^{r \times n_{\rm in}}$, $\bar{C} \in \real^{n_{\rm out} \times r}$,
where the dimension~$r$ is much smaller than~$n$.
Again we assume that $\bar{E}$ is non-singular.
Initial values $\bar{x}(0)=\bar{x}_0$ are approximated from
the initial values $x(0) = x_0$.
Nevertheless, the output of~(\ref{system-reduced}) should be a
good approximation to the output of~(\ref{linear-ode}),
i.e., $\bar{y}(t) \approx y(t)$ for all relevant times.
The system~(\ref{system-reduced}) is called the reduced-order model (ROM)
of the full-order model (FOM) given by~(\ref{linear-ode}).

The linear dynamical system~(\ref{system-reduced}) has its own
transfer function
$\bar{H} : \complex \backslash \bar{\Sigma} \rightarrow
\complex^{n_{\rm out} \times n_{\rm in}}$
of the form~(\ref{transfer}).
If both the original system~(\ref{linear-ode}) and the
reduced system~(\ref{system-reduced}) are asymptotically stable,
then error bounds are available in the case of $x_0=0$ and $\bar{x}_0=0$.
It holds that, see~\cite[p.~496]{benner-gugercin-willcox},
\begin{equation} \label{error-bound}
  \sup_{t \ge 0} \| y(t) - \bar{y}(t) \|_\infty \le
  \left\| H - \bar{H} \right\|_{\htwo} \| u \|_{\ltwo}
\end{equation} 
with the $\ltwo$-norm
\begin{equation} \label{l2-norm}
  \| u \|_{\ltwo} =
  \sqrt{ \int_0^{\infty} \| u(t) \|_2^2 \; {\rm d}t } \; , 
\end{equation}
the $\htwo$-norm~(\ref{h2-norm}),
the maximum (vector) norm $\| \cdot \|_\infty$ and the
Euclidean (vector) norm $\| \cdot \|_2$.

In projection-based MOR, see~\cite{antoulas},
two projection matrices $V,W \in \real^{n \times r}$ of full rank
are specified.
We obtain the matrices of the ROM~(\ref{system-reduced}) by
\begin{equation} \label{projected-matrices}
  \bar{A} = W^\top A V , \quad \bar{B} = W^\top B , \quad
  \bar{C} = C V \quad \bar{E} = W^\top E V .
\end{equation}
The orthogonality $V^\top V = I_r$ and sometimes the
biorthogonality $W^\top V = I_r$ are supposed
with the identity matrix $I_r \in \real^{r \times r}$.
Often the projection matrices result from the determination of
subspaces, i.e.,
\begin{equation} \label{subspaces}
  \mathcal{V} = {\rm span}(V) \subset \real^n
  \qquad \mbox{and} \qquad
  \mathcal{W} = {\rm span}(W) \subset \real^n .
\end{equation}
On the one hand, the original state variables are approximated within the
space~$\mathcal{V}$ by $x \approx V \bar{x}$.
On the other hand, the residual
\begin{equation} \label{residual}
  s(t) = E V \dot{\bar{x}}(t) - A V \bar{x}(t) - B u(t)
  \in \real^n
\end{equation}
is kept small by the requirement $s(t) \perp \mathcal{W}$
and thus $W^\top s(t) = 0$ for all~$t$.

It holds that $W=V$ in a Galerkin-type projection,
where just an appropriate projection matrix~$V$ has to be identified.
Important examples are the one-sided Arnoldi method, see~\cite{freund},
and proper orthogonal decomposition (POD), see~\cite{antoulas}.
However, the investigation of optimal choices for the matrix~$V$
is not within the scope of this paper.
We consider an arbitrary matrix~$V$ satisfying $V^\top V = I_r$.

In each MOR approach, the reduced system~(\ref{system-reduced}) is often
useless if it is not at least Lyapunov stable.
Many Krylov subspace methods or moment matching techniques
do not guarantee a stable ROM.


\subsection{Stability preservation by transformations}
\label{sec:preservation}
We apply the following well-known property for Galerkin-type
projection-based MOR, see~\cite{castane-selga,prajna}.

\begin{theorem} \label{thm:preservation}
  In the linear dynamical system~(\ref{linear-ode}),
  let $E$ be symmetric positive definite and $A$ be dissipative.
  It follows that the reduced system~(\ref{system-reduced}) is
  asymptotically stable for the matrices~(\ref{projected-matrices}) with a
  full-rank matrix~$V$ and $W=V$.
\end{theorem}

\underline{Proof:}

Let $B \equiv 0$ in~(\ref{linear-ode}) without loss of generality.
Since $E$ is symmetric positive definite and $V$ has full rank,
the matrix $\bar{E} = V^\top E V$ is also symmetric positive definite.
Let $\bar{E} = \bar{L} \bar{L}^\top$ be any symmetric decomposition,
for example, the Cholesky factorisation.
Thus the reduced system reads as
$\bar{L} \bar{L}^\top \dot{\bar{x}} = V^\top A V \bar{x}$,
which is equivalent to
\begin{equation} \label{z-dot}
  \dot{\bar{z}} = \bar{L}^{-1} V^\top A V \bar{L}^{-\top} \bar{z}
  = (V \bar{L}^{-\top})^\top A (V \bar{L}^{-\top}) \bar{z}
\end{equation}
with $\bar{z} = \bar{L}^\top \bar{x}$.
It follows that $V' = V \bar{L}^{-\top}$ is a full-rank matrix again.
We investigate the symmetric part of the transformed system matrix
in~(\ref{z-dot})
\begin{equation} \label{symm-part-A}
  V'^\top A V' + (V'^\top A V')^\top =
  V'^\top A V' + V'^\top A^\top V' = V'^\top (A + A^\top) V' .
\end{equation}
The dissipativity of~$A$ implies that $A + A^\top$ is negative definite
by Definition~\ref{def:diss-matrix}.
Due to the full rank of~$V'$,
the symmetric part~(\ref{symm-part-A}) is negative definite again.
The matrix $V'^\top A V'$ becomes dissipative and thus stable.
\hfill $\Box$

\medskip

An important special case of Theorem~\ref{thm:preservation}
is given for a system~(\ref{linear-ode}) with
an identity matrix $E=I_n$,
where explicit ODEs arise.

We assume just a non-singular matrix~$E$ now.
If $E$ is symmetric and positive definite, then we suppose
a non-dissipative matrix~$A$.
Hence Theorem~\ref{thm:preservation} is not applicable.
The idea is to transform the system~(\ref{linear-ode}) into
an equivalent system with a symmetric positive definite mass matrix
and a dissipative matrix on the right-hand side.
Let $M \in \real^{n \times n}$ be symmetric and positive definite.
We multiply the system~(\ref{linear-ode}) by the non-singular matrix
$E^\top M$ and obtain the equivalent system
\begin{equation} \label{ode-trafo}
  E^\top M E \dot{x}(t) = E^\top M A \dot{x} + E^\top M B u(t) .
\end{equation}
This operation can be seen as a transformation in the range of
the linear mappings.
The mass matrix $E^\top M E$ is always symmetric and positive definite.
The matrix~$M$ has to be chosen such that the matrix
$E^\top M A$ is dissipative.
This task is achieved by solving the generalised Lyapunov equations,
cf.~\cite[p.~449]{mueller},
\begin{equation} \label{lyapunov}
  A^\top M E + E^\top M A + F = 0
\end{equation}
for any symmetric positive definite matrix $F \in \real^{n \times n}$.
Since $E^{-1}A$ is a stable matrix,
the Lyapunov equations~(\ref{lyapunov}) exhibit a unique solution~$M$,
which is also symmetric positive definite.
Moreover, the matrix $\hat{A}=E^\top M A$ is dissipative
in the system~(\ref{ode-trafo}),
because $\hat{A} + \hat{A}^\top = -F$ is negative definite.
We summarise this important result in a theorem.

\begin{theorem} \label{thm:stable}
  Let the dynamical system~(\ref{linear-ode}) be asymptotically stable.
  If $M$ is the solution of the Lyapunov equations~(\ref{lyapunov})
  including a symmetric positive definite matrix~$F$,
  then each Galerkin-type projection-based MOR of
  the linear dynamical system~(\ref{ode-trafo}) 
  yields an asymptotically stable reduced system~(\ref{system-reduced}).
\end{theorem}

The proof follows from Theorem~\ref{thm:preservation} and
the discussion above.

\subsection{Numerical solution of Lyapunov equations}
There are two classes of numerical linear algebra methods for the solution of
the (generalised) Lyapunov equation~(\ref{lyapunov}):
\begin{enumerate} 
\item[i)]
  Direct methods:
  Either the solution~$M$ or its Cholesky factor $L$ without solving for~$M$
  first can be computed, see~\cite{hammarling}.
  The computational effort reads as $\mathcal{O}(n^3)$.
  The memory requirement is about $\frac{n^2}{2}$ machine numbers.
\item[ii)]
  Approximate methods of the following types, see~\cite{kramer,wolf}:
  \begin{itemize}
  \item
    projection techniques: Krylov subspace methods, POD, etc.,
  \item
    alternating direction implicit (ADI) method.
  \end{itemize}
  The approximate methods yield a low-rank factor
  $Z \in \real^{n \times q}$ with $q \ll n$ satisfying $M \approx Z Z^\top$.
\end{enumerate}
On the one hand, the dimension~$n$ is high in MOR and thus the
direct methods are often computationally infeasible or their
effort is much higher than solving the FOM~(\ref{linear-ode}).
On the other hand, approximate techniques typically require much less
computational effort.
In addition, the decomposition $M \approx ZZ^\top$ with a low-rank factor~$Z$
allows for cheap matrix-matrix multiplications with~$M$ in the projections. 
However, we encounter two difficulties within this approach:
\begin{enumerate}
\item
  The approximation $\tilde{M} = ZZ^\top$ is singular by construction
  and thus does not describe a basis transformation.
\item
  The Lyapunov equations~(\ref{lyapunov}) require a positive definite
  matrix~$F$ due to Theorem~\ref{thm:preservation}.
  Yet the ADI iteration is efficient only in the case of
  $F = Z_F Z_F^\top$ with a low-rank factor $Z_F \in \real^{n \times q_F}$
  satisfying $q_F \ll n$, see~\cite[p.~10]{penzl}.
  Now $F$ is just positive semi-definite.
\end{enumerate}  
We discuss the first item further.
The mass matrix reads as
\begin{equation} \label{reduced-E}
  \bar{E} = V^\top E^\top \tilde{M} E V = (Z^\top EV)^\top (Z^\top EV)
\end{equation}
in the reduced system~(\ref{system-reduced}).
Since $V$ has full rank, $EV \in \real^{n \times r}$
is also a full-rank matrix. 
The subspace
\begin{equation} \label{range}
  \mathcal{R}(S) = \{ Sx \, : \, x \in \real^m \} \subseteq \real^n
  \qquad \mbox{for} \;\; S \in \real^{n \times m}
\end{equation}
denotes the range of the linear mapping induced by a matrix.
If it holds that
$\mathcal{R}(Z) \perp \mathcal{R}(EV)$,
then we obtain $Z^\top E V = 0$ and thus $\bar{E} = 0$.
Although this event may happen due to $q,r \ll n$,
it is rather unlikely.
Nevertheless, the matrix~(\ref{reduced-E}) can easily become ill-conditioned
or singular with respect to working precision.

We observe that the stabilisation is not straightforward.
An alternative technique is suggested to omit the two difficulties above
in the next section.


\section{Numerical method}
\label{sec:method}
A numerical approach is derived and investigated, which guarantees
the preservation of stability in a Galerkin-type projection-based MOR,
while Lyapunov equations can be solved approximately.

\subsection{Setup of projection matrices}
\label{sec:preliminaries}
Let $V \in \real^{n \times r}$ with $V^\top V = I_r$ be a projection matrix
constructed for a reduction of the system~(\ref{linear-ode}).
We apply the Galerkin-type MOR with~$V$ to the
transformed system~(\ref{ode-trafo}).
The matrices of the associated reduced system~(\ref{system-reduced})
can be written in the form~(\ref{projected-matrices})
with the projection matrix
\begin{equation} \label{matrix-W}
  W = M E V .
\end{equation}
Hence this reduction is of the form~(\ref{projected-matrices})
and thus represents a special case of a (non-Galerkin type)
MOR for the original system~(\ref{linear-ode}).
However, it holds that $W^\top V \neq I_r$ in general, i.e.,
biorthogonality is not given.

We briefly discuss the computation work to obtain the 
matrices~(\ref{projected-matrices}).
The mass matrix~$E$ is often more sparse than the matrix~$A$.
Thus the matrix-matrix multiplication $EV$ is relatively cheap.
The calculation of~(\ref{matrix-W}) requires $W = M(EV)$,
where $M$ is dense.
We will use an approximation of~$M$ in Section~\ref{sec:lyapunov},
which allows for a cheap computation of this matrix-matrix product.
Given $V$ and $W$, the matrices~(\ref{projected-matrices})
can be computed as in any projection-based MOR.

It is important to note that the two Galerkin-type MORs with projection
matrix~$V$ for the system~(\ref{linear-ode}) and the transformed
system~(\ref{ode-trafo}) are not equivalent.
The reason is that the residual~(\ref{residual}) is mapped to
$$ E^\top M  s(t) =
E^\top M E V \dot{\bar{x}}(t) - E^\top M A V \bar{x}(t) - E^\top M B u(t)
\qquad \mbox{for} \;\; t \ge 0 . $$
The property $s(t) \perp \mathcal{V}$ is not equivalent to
the property $E^\top M s(t) \perp \mathcal{V}$ for each~$t$.

\subsection{Alternative ansatz and Lyapunov equations}
\label{sec:lyapunov}
In view of Theorem~\ref{thm:stable}, the construction of the
projection matrices involves the efficient solution of the
Lyapunov equation~(\ref{lyapunov}) together with the identification
of an appropriate matrix~$F$.
Concerning the system~(\ref{linear-ode}), 
we assume that the matrix $E^{-1}A$ is stable but non-dissipative.
Definition~\ref{def:diss-matrix} implies that the symmetric part
\begin{equation} \label{symm-part}
  G_{\rm sym} = E^{-1}A + A^\top E^{-\top}
\end{equation}
has $k \ge 1$ non-negative eigenvalues.
Since the matrix~$E^{-1}A$ has exclusively eigenvalues with negative real part,
we expect a small number $k$.
Let $\mu_1 \ge \mu_2 \ge \cdots \ge \mu_n$ be the eigenvalues
of~(\ref{symm-part}) and $u_1,u_2,\ldots,u_n$ the associated eigenvectors,
which form an orthonormal basis.
It holds that $\mu_{\max} = \mu_1$ is the maximum eigenvalue.
We define $U = (u_1,\ldots,u_k) \in \real^{n \times k}$.
The following lemma motivates a choice of the matrix~$F$
in~(\ref{lyapunov}).

\begin{lemma} \label{lemma:posdef}
  Let $A,E \in \real^{n \times n}$, where $E$ is non-singular
  and $E^{-1} A$ is a non-dissipative matrix.
  It follows that the symmetric matrix
  \begin{equation} \label{matrix-F}
    F = - (E^{-1}A + A^\top E^{-\top}) + (\mu_{\max} + \delta) U U^\top
  \end{equation}
  with any real number $\delta > 0$ is positive definite.
\end{lemma}

\underline{Proof:} \nopagebreak

The eigenvectors $u_1,\ldots,u_n$ of the symmetric matrix~(\ref{symm-part})
represent an orthonormal basis.
It follows that
$$ F u_j = - G_{\rm sym} u_j + (\mu_1+\delta) U U^\top u_j =
- \mu_j u_j + (\mu_1 + \delta) u_j = (\mu_1 - \mu_j + \delta) u_j $$
for $j=1,\ldots,k$
and
$$ F u_j = - G_{\rm sym} u_j + (\mu_1+\delta) U U^\top u_j = - \mu_j u_j $$
for $j=k+1,\ldots,n$.
Hence the eigenvectors of~$F$ are $u_1,\ldots,u_n$ again.
Due to $\mu_1 \ge \mu_j$ for all~$j$ and $\mu_j<0$ for $j>k$
and $\mu_1 \ge 0$,
all eigenvalues of~$F$ are positive.
\hfill $\Box$

\medskip

The matrix~(\ref{matrix-F}) represents a rank-$k$ update of
a matrix, where the eigenvalue problem was investigated
in~\cite{ding-yao}, for example.

Concerning the solution of the Lyapunov equation~(\ref{lyapunov}),
we perform the ansatz
\begin{equation} \label{ansatz-M}
  M = E^{-\top} E^{-1} + \Delta M . 
\end{equation}
The simplification $M = I_n + \Delta M$ arises in the case of $E = I_n$.
Furthermore, an approximation of~(\ref{ansatz-M}) reads as
\begin{equation} \label{appr-M}
  \tilde{M} = E^{-\top} E^{-1} + Z Z^\top
\end{equation}
including a low-rank factorisation with a factor
$Z \in \real^{n \times q}$ of rank~$q < n$.
Hence the approximation consists in $\Delta M \approx Z Z^\top$.

Let
\begin{equation} \label{ansatz-F}
  F = -( E^{-1} A + A^\top E^{-\top}) + \Delta F .
\end{equation}
We require a positive definite matrix~$F$.
Inserting~(\ref{ansatz-M}) and~(\ref{ansatz-F})
in~(\ref{lyapunov}),
we obtain the alternative Lyapunov equation
\begin{equation} \label{lyapunov2}
  A^\top \Delta M E + E^\top \Delta M A + \Delta F = 0 .
\end{equation}
Due to Lemma~\ref{lemma:posdef}, we choose
\begin{equation} \label{delta-F}
  \Delta F = (\mu_{\max}+\delta) U U^\top ,
\end{equation}
which guarantees that the matrix~(\ref{ansatz-F}) is positive definite.
The Lyapunov equations~(\ref{lyapunov2}) can be written as
\begin{equation} \label{lyapunov3}
  A^\top \Delta M E + E^\top \Delta M A + \tilde{U} \tilde{U}^\top = 0 
\end{equation}
with the factor
$\tilde{U} = \sqrt{\mu_{\max}+\delta} \, U \in \real^{n \times k}$.
An efficient approximate solution of the Lyapunov
equations~(\ref{lyapunov3}) with factorisation $\Delta M \approx Z Z^\top$
and $Z \in \real^{n \times q}$ ($q \ge k$)
is feasible in the case of low numbers~$k$.
The determination of the matrix~$\tilde{U}$ requires the computation
of the eigenvectors associated to the~$k$ dominant eigenvalues 
for a real symmetric matrix.
Appropriate iterative methods exist for the numerical computation
of just the dominant eigenvalues and their eigenvectors.

The ROM~(\ref{system-reduced}) includes the mass matrix~$\bar{E}$.
We obtain a bound on its condition number.


\begin{theorem} \label{thm:condition}
  The matrix $\bar{E} = V^\top E^\top \tilde{M} E V$
  including the approximation~(\ref{appr-M}) and $V^\top V = I_r$
  is symmetric positive definite. 
  The condition number with respect to the spectral (matrix) norm
  $\| \cdot \|_2$ exhibits the upper bound
  \begin{equation} \label{cond-bound}
    {\rm cond} (\bar{E}) \le 1 + \| E \|_2^2 \| Z \|_2^2 .
  \end{equation}
\end{theorem}

\underline{Proof:}

If holds that
$$ \bar{E} = V^\top E^\top (E^{-\top} E^{-1} + Z Z^\top) E V =
I_r + (V^\top E^\top Z) (V^\top E^\top Z)^\top . $$
Obviously, the matrix~$\bar{E}$ is symmetric.
Let $y \in \real^r$ be an eigenvector of $\bar{E}$ with $\| y \|_2=1$
and $\eta$ be the associated eigenvalue.
It follows that
$$ \eta = y^\top \bar{E} y =
1 + \| Z^\top E V y \|_2^2 =
\| y \|_2^2 + \| Z^\top E V y \|_2^2 \le
1 + \| E \|_2^2 \| Z \|_2^2 $$
due to $\| V \|_2 = 1$ and $\| Z^\top \|_2 = \| Z \|_2$.
Thus the eigenvalues are bounded uniformly by
$1 \le \eta \le 1+\| E \|_2^2 \| Z \|_2^2$.
All eigenvalues are positive.
Since the matrix is symmetric, the condition number satisfies
$$ {\rm cond} (\bar{E}) = \frac{\eta_{\max}}{\eta_{\min}} \le
1+\| E \|_2^2 \| Z \|_2^2 , $$
which completes the proof.
\hfill $\Box$

\medskip

Hence the upper bound~(\ref{cond-bound}) depends only on the magnitude
of the factors~$E$, $Z$ and not on the position of the
subspace $\mathcal{R}(Z)$ in~$\real^n$.

\subsection{Symmetric decomposition}
\label{sec:decomposition}
In~\cite{prajna}, a symmetric decomposition of the solution of
Lyapunov equations is used to perform a basis transformation.
In Section~\ref{sec:lyapunov}, we introduced an approach,
which does not require a symmetric decomposition of the
approximation~(\ref{appr-M}).
Yet a numerical technique is feasible including such a
symmetric decomposition with a reasonable computational effort
in the case of $E=I_n$ for~(\ref{linear-ode}).
We outline this scheme, even though it is not used for the numerical
computations in Section~\ref{sec:example}.

It holds that $\tilde{M} x \in \mathcal{R}(Z)$ for
$x \in \mathcal{R}(Z)$ and
$\tilde{M} x = x$ for $x \in \mathcal{R}(Z)^\perp$
concerning subspaces~(\ref{range}).
Due to $q \ll n$, the linear mapping described by the matrix~$\tilde{M}$
represents the identity on a relatively large subspace.
We require a partition
$\real^n = \mathcal{R}(Z) \oplus \mathcal{R}(Z)^\perp$.
Therefore, we apply the matrix square root for the symmetric decomposition
$\tilde{M} = \tilde{M}^{\frac{1}{2}} \tilde{M}^{\frac{1}{2}}$
of~(\ref{appr-M}).
Note that the Cholesky algorithm is inappropriate, because it 
calculates the complete matrix $Z Z^\top$ in~(\ref{appr-M})
followed by a computational effort of $\mathcal{O} (n^3)$.

\begin{theorem} \label{thm:matrix-root}
  Consider the $QR$-factorisation
  \begin{equation} \label{qr-decomp}
    Z = Q R \quad \mbox{with} \quad
    R = \begin{pmatrix} R' \\ 0 \\ \end{pmatrix} ,
  \end{equation}
  where $Q \in \real^{n \times n}$ is an orthogonal matrix,
  $R \in \real^{n \times q}$
  and $R' \in \real^{q \times q}$ is an upper triangular matrix.
  The eigendecomposition
  \begin{equation} \label{eigen-R}
    R' R'^\top = S D S^\top
  \end{equation}
  includes an orthogonal matrix~$S \in \real^{q \times q}$ and a
  diagonal matrix~$D \in \real^{q \times q}$ with positive diagonal elements.
  The matrix square root of~(\ref{appr-M}) exhibits the formula
  \begin{equation} \label{M-root}
    \tilde{M}^\frac{1}{2} =
    Q \tilde{S} \tilde{D}^\frac{1}{2} \tilde{S}^\top Q^\top
  \end{equation}
  with the matrices
  $$ \tilde{S} = \begin{pmatrix} S & 0 \\ 0 & I_{n-q} \\ \end{pmatrix}
  \qquad \mbox{and} \qquad
  \tilde{D} =
  \begin{pmatrix} (I_q+D)^\frac{1}{2} & 0 \\ 0 & I_{n-q} \\ \end{pmatrix} $$
  including identity matrices $I_j \in \real^{j \times j}$
  for $j=q,n-q$.
\end{theorem}

\underline{Proof:} \nopagebreak 

Since we always assume ${\rm rank}(Z)=q$,
it follows that ${\rm rank}(R)={\rm rank}(R')=q$.
Hence the diagonal elements of~$D$ are strictly positive.
We calculate
$$ \begin{array}{rcl}
  \tilde{M} & = & I_n + Z Z^\top \;\; = \;\; I_n + (QR)(QR)^\top \\[2ex]
  & = & I_n + Q \begin{pmatrix} R' R'^\top & 0 \\ 0 & 0 \\ \end{pmatrix} Q^\top
  \; = \;
  I_n + Q \begin{pmatrix} S D S^\top & 0 \\ 0 & 0 \\ \end{pmatrix} Q^\top \\[3ex]
  & = &
  I_n + Q \tilde{S} \begin{pmatrix} D & 0 \\ 0 & 0 \\ \end{pmatrix}
  \tilde{S}^\top Q^\top
  \; = \;
  (Q \tilde{S})
  \begin{pmatrix} I_q + D & 0 \\ 0 & I_{n-q} \\ \end{pmatrix}
  (Q \tilde{S})^\top . \\
  \end{array} $$
We have achieved an orthogonal eigendecomposition of~$\tilde{M}$.
Taking the square roots of the elements in the diagonal matrix yields
the desired matrix square root~(\ref{M-root}).
\hfill $\Box$

\medskip

The matrix square root is a symmetric matrix again.
However, the basis change induced by $\tilde{M}^\frac{1}{2}$
is not an orthogonal transformation.
We directly obtain the inverse matrix of~(\ref{M-root}) by
$$ \tilde{M}^{-\frac{1}{2}} =
Q \tilde{S} \tilde{D}^{-\frac{1}{2}} \tilde{S}^\top Q^\top , $$
where just reciprocal values are required within the diagonal matrix.

Householder transformations perform the $QR$-factorisation of
a matrix, see~\cite[p.~223]{stoerbulirsch}.
Thus the matrix~$Q$ is given by successive rank-one updates.
The complexity becomes $2nq^2 - \frac{2}{3}q^3$ operations
for~(\ref{qr-decomp}), which can be characterised as $\mathcal{O}(n q^2)$
due to $q \ll n$.
Likewise, the eigendecomposition~(\ref{eigen-R}) is cheap,
because matrices of dimension~$q$ are involved.

Basis transformations of projection matrices~$V$ require
matrix-matrix multiplications, which can be performed by just
a few matrix-vector products
$\tilde{M}^{\frac{1}{2}} v$ and $\tilde{M}^{-\frac{1}{2}} v$.
We consider $\tilde{M}^{\frac{1}{2}} v$ without loss of generality.
The computational effort consists of a sequence of matrix-vector
multiplications:
\begin{enumerate}
\item
  $v^{(1)} = Q^\top v$:
  The matrix~$Q$ is given by $q$~Householder transformations.
\item
  $v^{(2)} = \tilde{S}^\top v^{(1)}$:
  The computation work mainly is a matrix-vector
  multiplication with the dense matrix~$S^\top$ of dimension~$q$
  and thus $\mathcal{O}(q^2)$.
\item
  $v^{(3)} = \tilde{D}^\frac{1}{2} v^{(2)}$: The computational effort is only
  $\mathcal{O}(q)$.
\item
  $v^{(4)} = \tilde{S} v^{(3)}$: As in step~2.
\item
  $v^{(5)} = Q v^{(4)}$: As in step~1.
\end{enumerate}
In conclusion, the effort for each matrix-vector product is negligible
in comparison to the combination of $QR$-factorisation~(\ref{qr-decomp})
and eigendecomposition~(\ref{eigen-R}).


\section{Application to nonlinear dynamical systems}
\label{sec:nonlinear}
We show the applicability of the stability-preserving technique
from Section~\ref{sec:method} to nonlinear systems of ODEs.

\subsection{Projection-based model order reduction}
Let a nonlinear system of ODEs be given in the form
\begin{equation} \label{nonlinear-system}
  \begin{array}{rcl}
    E \dot{x}(t) & = & f(x(t)) + B u(t) \\
    y(t) & = & C x(t) \\
  \end{array}
\end{equation}
with state variables $x : [0,\infty) \rightarrow \real^n$,
inputs $u : [0,\infty) \rightarrow \real^{n_{\rm in}}$, 
outputs $y : [0,\infty) \rightarrow \real^{n_{\rm out}}$,
matrices $B \in \real^{n \times n_{\rm in}}$ and $C \in \real^{n_{\rm out} \times n}$, 
a non-singular mass matrix $E \in \real^{n \times n}$ and a nonlinear
smooth function $f : \real^n \rightarrow \real^n$.
Again initial values $x(0) = x_0$ are predetermined.
Concerning general nonlinear dynamical system,
the input-output behaviour cannot be described by a transfer function
in the frequency domain.

A projection-based MOR involves matrices $V,W \in \real^{n \times r}$
of full rank again.
The ROM reads as
\begin{equation} \label{rom-nonlinear}
  \begin{array}{rcl}
    \bar{E} \dot{\bar{x}}(t) & = & \bar{f}(\bar{x}(t)) + \bar{B} u(t) \\
    \bar{y}(t) & = & \bar{C} \bar{x}(t) . \\
  \end{array}
\end{equation}
The matrices are defined by~(\ref{projected-matrices}).
The nonlinear function is approximated by
\begin{equation} \label{reduced-f}
  \bar{f}(\bar{x}) = W^\top f( V \bar{x} )
  \qquad \mbox{for} \;\; \bar{x} \in \real^r .
\end{equation}
Again the choice $W=V$ yields a method of Galerkin-type.

\subsection{Stability of stationary solutions}
We assume an autonomous system of ODEs~(\ref{nonlinear-system})
with a constant input $u(t) \equiv u_0$.
Without loss of generality, we drop the term $Bu$ in~(\ref{nonlinear-system}),
because $Bu_0$ can be included in the function~$f$.
A stationary solution or equilibrium $x^* \in \real^n$
of the nonlinear dynamical system is characterised by $f(x^*) = 0$.
We assume isolated stationary solutions.
The stationary solution is asymptotically stable, if and only if the matrix
\begin{equation} \label{jacobian} \textstyle
  A' = E^{-1} \left. \frac{\partial f}{\partial x} \right|_{x=x^*}
\end{equation}
is a stable matrix with respect to Definition~\ref{def:stable-matrix},
see~\cite[p.~22]{seydel}.
We assume an equilibrium $x^* = 0 \in \real^n$ without loss of generality.
It holds that $f(0)=0$.

In a projection-based MOR, the reduced system~(\ref{rom-nonlinear})
features the stationary solution $\bar{x}^* = 0 \in \real^r$ now.
The stability of this equilibrium is determined by the matrix
$$ \bar{A}' = \bar{E}^{-1} \textstyle
\left. \frac{\partial \bar{f}}{\partial \bar{x}} \right|_{\bar{x}=0} =
\bar{E}^{-1} W^\top
\left. \frac{\partial f}{\partial x} \right|_{x=0} V $$
due to (\ref{reduced-f}) and (\ref{jacobian}).
Yet the stationary solution may be unstable, even if the
zero-state of the FOM~(\ref{nonlinear-system}) is asymptotically stable.

\subsection{Preservation of stability}
We mimic the method from Section~\ref{sec:preservation}.
Let $M \in \real^{n \times n}$ be symmetric and positive definite.
A multiplication of the system~(\ref{nonlinear-system}) by $E^\top M$
yields (assuming $Bu=0$)
\begin{equation} \label{nonlinear-trafo}
  E^\top M E \dot{x}(t) = E^\top M f(x(t)) .
\end{equation}
The state $x^*=0$ also represents a stationary solution
of~(\ref{nonlinear-trafo}).
We arrange the Lyapunov equation
\begin{equation} \label{lyapunov-nonlinear}
  \textstyle
  \left( \left. \frac{\partial f}{\partial x} \right|_{x=0} \right)^\top M E
  + E^\top M \left( \left. \frac{\partial f}{\partial x} \right|_{x=0} \right)
  + F = 0
\end{equation}
including a symmetric positive definite matrix~$F \in \real^{n \times n}$.
The following theorem summarises the stability-preserving approach.

\begin{theorem} \label{thm:nonlinear}
  Let the nonlinear dynamical system~(\ref{nonlinear-system})
  with $Bu=0$ have the asymptotically stable stationary solution $x^*=0$.
  Let~$V \in \real^{n \times r}$ be any projection matrix with full rank. 
  The Lyapunov equation~(\ref{lyapunov-nonlinear}) yields a solution~$M$
  for some symmetric positive definite matrix~$F$.
  It follows that the reduced system~(\ref{rom-nonlinear}) given
  by~(\ref{projected-matrices}) and~(\ref{reduced-f}) with the
  matrix~$W$ from~(\ref{matrix-W})
  features a stationary solution~$\bar{x}^*=0$,
  which is asymptotically stable.
\end{theorem}

\underline{Proof:}

The projections~(\ref{projected-matrices}) and~(\ref{reduced-f})
with the matrix~(\ref{matrix-W}) generate the reduced system
$$ V^\top E^\top M E V \dot{\bar{x}}(t)
= V^\top E^\top M f (V \bar{x}(t)) . $$
The matrix $\bar{E} = V^\top E^\top M E V$ is symmetric and positive definite.
Let $\bar{E} = \bar{L} \bar{L}^\top$.
Using $\bar{z} = \bar{L}^\top \bar{x}$, we obtain the equivalent
dynamical system
$$ \dot{\bar{z}}(t) =
( V \bar{L}^{-\top} )^\top E^\top M f ( V \bar{L}^{-\top} \bar{z}(t) ) . $$
The stability of the stationary solution $\bar{z}^* = \bar{x}^* = 0$
is determined by the matrix
$$ \bar{A}'' =  ( V \bar{L}^{-\top} )^\top E^\top M 
\left. \textstyle \frac{\partial f}{\partial x} \right|_{x = 0} 
( V \bar{L}^{-\top} ) . $$
Since $M$ is the solution of the Lyapunov equation~(\ref{lyapunov-nonlinear}),
it follows that
$E^\top M \left. \frac{\partial f}{\partial x} \right|_{x=0}$
is dissipative.
As in the proof of Theorem~\ref{thm:preservation},
we conclude that $\bar{A}''$ is dissipative and thus a stable matrix.
\hfill $\Box$

\medskip

Now the ansatz~(\ref{ansatz-M}) and the approximation~(\ref{appr-M}) can
be applied to solve the Lyapunov equation~(\ref{lyapunov-nonlinear})
approximately.
The technique proceeds as in Section~\ref{sec:lyapunov}.
In contrast to~\cite{prajna}, a transformation is not required in the
state space and thus a symmetric decomposition of the solution
to the Lyapunov equation~(\ref{lyapunov-nonlinear}) can be omitted.


\section{Illustrative examples}
\label{sec:example}
We apply the approach from Section~\ref{sec:problem-def} and
Section~\ref{sec:method} to two examples now.
All numerical computations were performed by the software package
MATLAB~\cite{matlab2016}.

\subsection{Stochastic model of mass-spring-damper system}
Lohmann and Eid~\cite{lohmann-eid} introduced a mass-spring-damper
configuration, which can be modelled by a first-order
system~(\ref{linear-ode}) with $n'=8$ equations.
The system is single-input-single-output (SISO).
In~\cite{pulch18}, this example was extended to a stochastic model,
where physical parameters are replaced by random variables 
with uniform distributions.
A polynomial chaos expansion, see~\cite{xiu-book}, is used
with $m$~orthogonal basis polynomials in the random space.
A stochastic Galerkin method yields a linear dynamical
system~(\ref{linear-ode}) of dimension $n=m n'$.
The system is single-input-multiple-output (SIMO).
More details can be found in~\cite{pulch18}.

\begin{figure}
  \begin{center}
  \includegraphics[width=6.5cm]{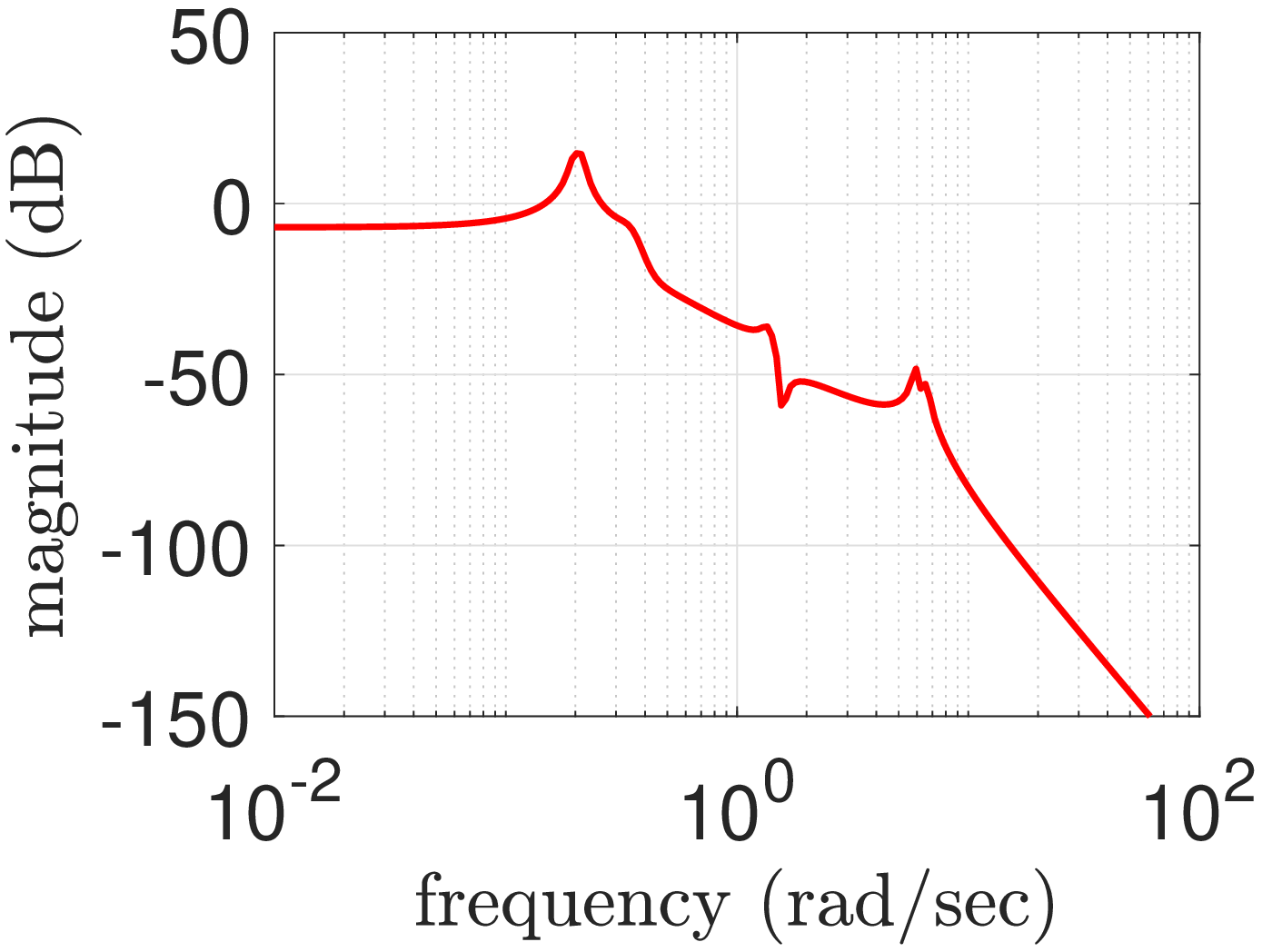}
  \hspace{5mm}
  \includegraphics[width=6.5cm]{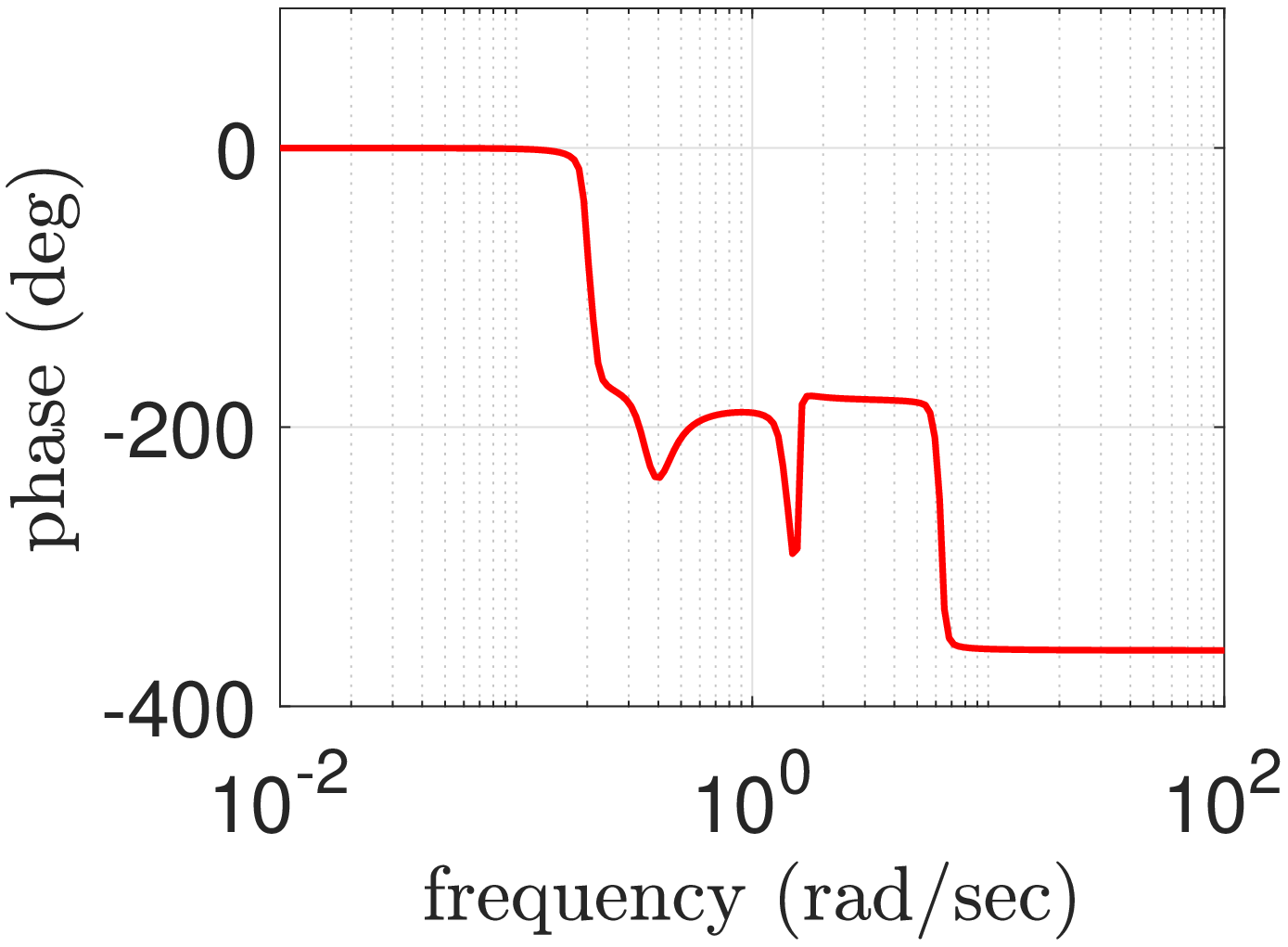}
  \end{center}
  \caption{Bode plot of stochastic Galerkin system for mass-spring-damper
    configuration.}
\label{fig:msd-bode}
\end{figure}

We consider this example again, where the uniformly distributed random
parameters exhibit a variation of 10\% around their mean values.
All basis polynomials up to degree three are included ($m=680$).
In the system~(\ref{linear-ode}),
the mass matrix~$E$ is symmetric positive definite and
the matrix $E^{-1}A$ is stable.
However, $A$ is not dissipative and thus Theorem~\ref{thm:preservation}
cannot be applied.
Table~\ref{tab:msd} shows further properties of the
linear dynamical system.
We focus on the constant basis polynomial, which is associated to
an approximation of the expected value of the original single output.
Thus our stochastic Galerkin system becomes SISO,
whose Bode plot is depicted in Figure~\ref{fig:msd-bode}.

\begin{table}
  \caption{Properties of stochastic mass-spring-damper system.\label{tab:msd}}
\begin{center}
  \begin{tabular}{cc} \hline
    dimension~$n$ & 5440 \\
    \# non-zero entries in~$A$ & 25120 \\
    \# non-zero entries in~$E$ & 6400 \\ \hline
  \end{tabular}
\end{center}
\end{table}

We want to employ a transformation to guarantee the stability in ROMs.
Yet the symmetric part~(\ref{symm-part}) of the matrix $E^{-1}A$ exhibits
$k=2720$ non-negative eigenvalues,
which violates the assumption of $k \ll n$.
Thus we do not use the approximation~(\ref{appr-M}).
Alternatively, we directly solve the generalised Lyapunov
equation~(\ref{lyapunov}) by the MATLAB function {\tt lyap},
where the identity matrix~$F=I_n$ is chosen.
Given any projection matrix~$V$, the stability-preserving reduction
is done via~(\ref{projected-matrices}),(\ref{matrix-W}).

\begin{figure}
  \begin{center}
  \includegraphics[width=6.5cm]{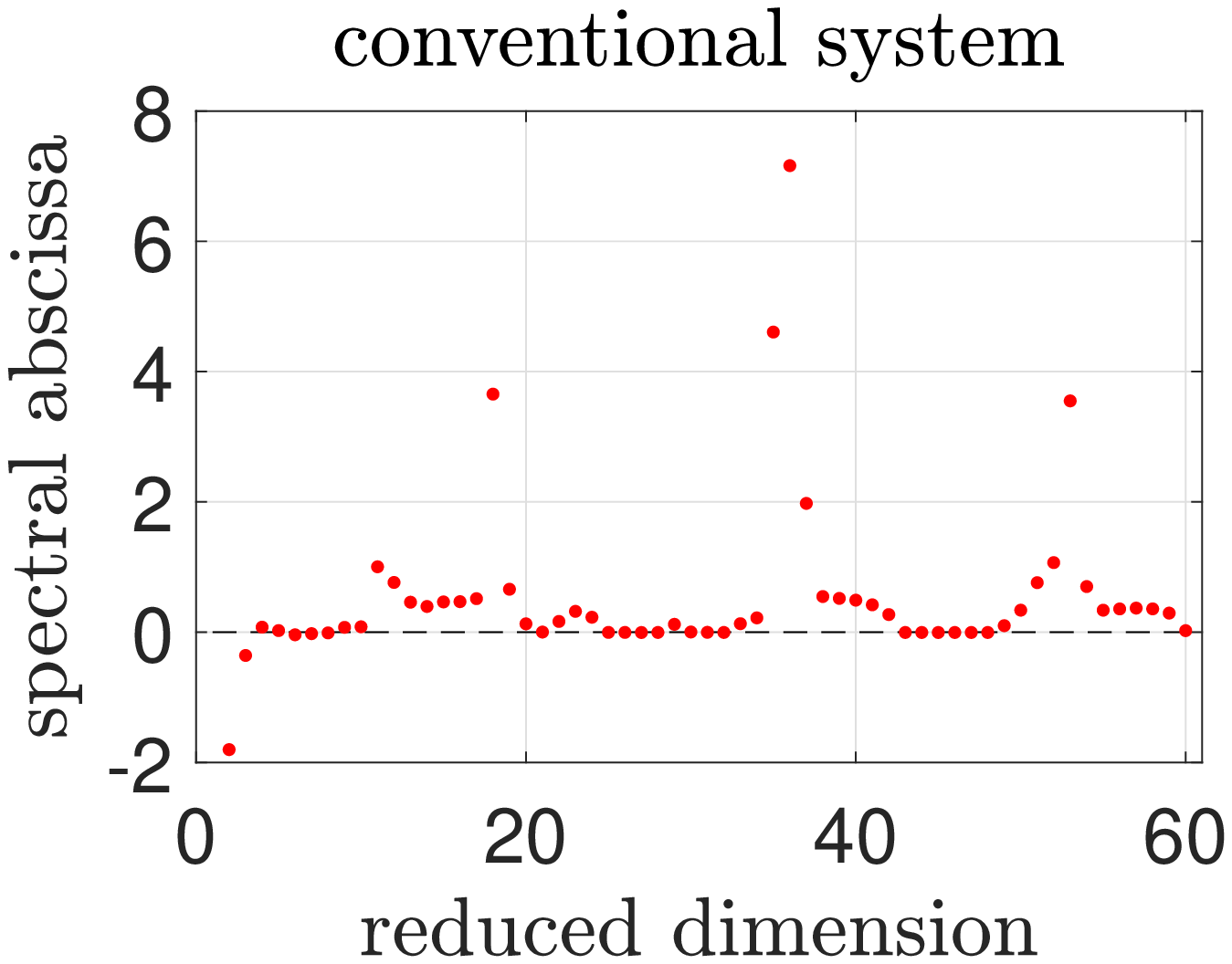}
  \hspace{5mm}
  \includegraphics[width=6.5cm]{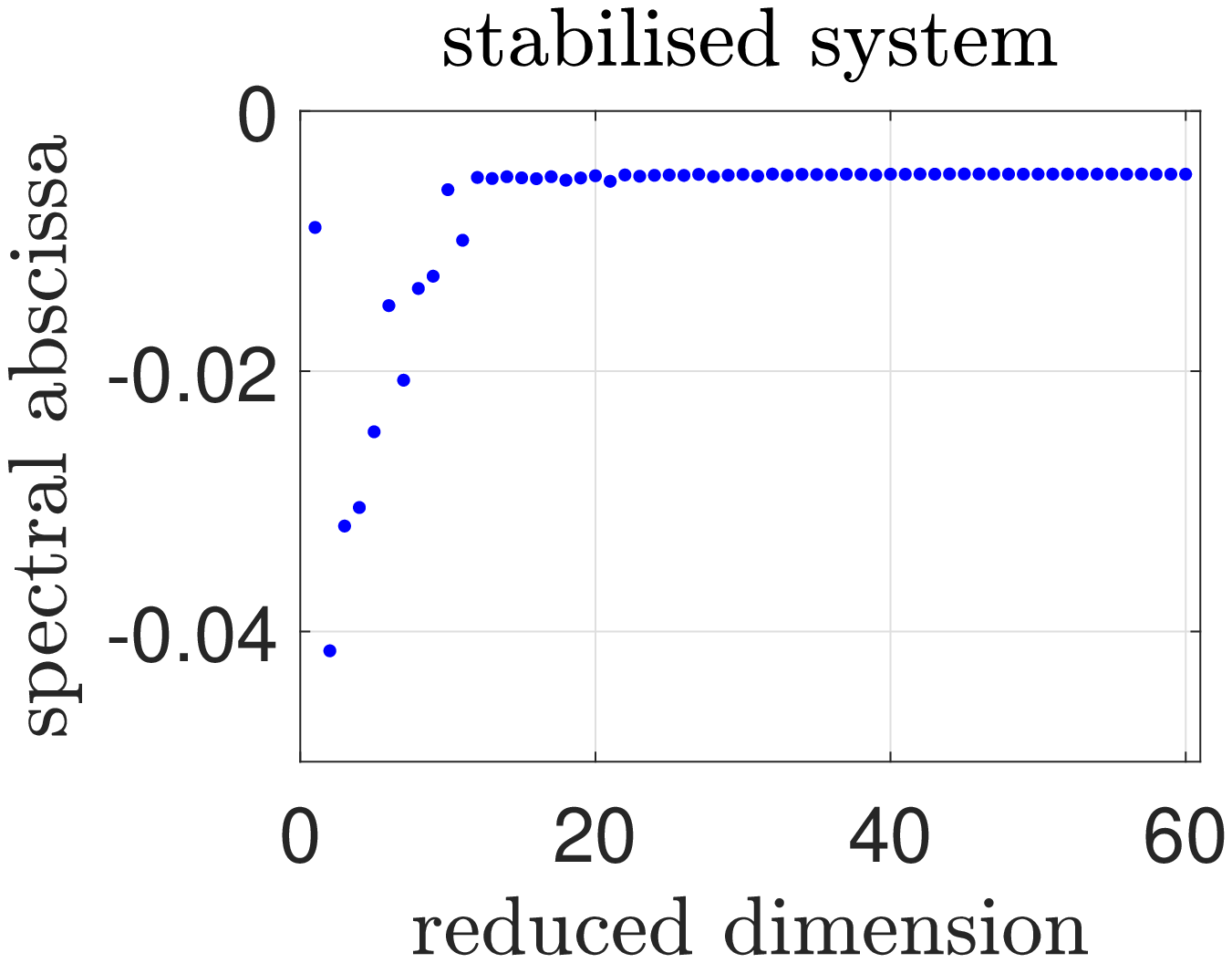}
  \end{center}
  \caption{Spectral abscissa of the matrices in the ROMs
    from conventional system (left)
    and stabilised system (right).}
\label{fig:msd-spectrum}
\end{figure}

The one-sided Arnoldi algorithm with the real expansion point $s_0=1$
yields the projection matrices $V \in \real^{n \times r}$ for each
integer~$r$ now.
We examine the cases $r=1,2,\ldots,60$.
On the one hand, the ROM is determined by the conventional
form~(\ref{projected-matrices}) with $W=V$.
Figure~\ref{fig:msd-spectrum} (left) illustrates the spectral abscissa
from Definition~\ref{def:spectral} for the matrices $\bar{E}^{-1} \bar{A}$.
It follows that just 18 of the 60 reduced systems~(\ref{system-reduced})
are stable.
On the other hand, the stability-preserving reduction
(\ref{projected-matrices}) with~(\ref{matrix-W}) is arranged.
Figure~\ref{fig:msd-spectrum} (right) shows the associated spectral abscissas.
As expected, all reduced systems~(\ref{system-reduced}) are stable
in this method.


\begin{figure}
\begin{center}
\includegraphics[width=8cm]{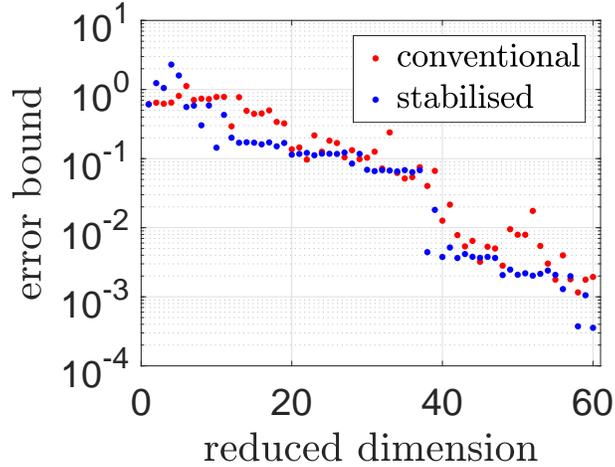} 
\end{center}
\caption{Error bound in $\htwo$-norm for the two MOR approaches
  in mass-spring-damper example.}
\label{fig:msd-error}
\end{figure}

We also compare the approximation quality between the
conventional reduced systems and the stabilised reduced systems,
because the two approaches are not equivalent even if the systems are stable.
Therefore we compute the error bounds~(\ref{error-bound})
for a unit norm~(\ref{l2-norm}) of the input ($\| u \|_{\ltwo} = 1$),
which are the $\htwo$-norms~(\ref{h2-norm}) of the differences between the
transfer functions of original system and reduced system.
Figure~\ref{fig:msd-error} depicts approximations of these upper bounds
computed using a grid on the imaginary axis.
We observe that the magnitudes of the error indicators are often
lower and thus better in the stabilised systems.

\subsection{Anemometer benchmark}
\begin{figure}
  \begin{center}
  \includegraphics[width=6.5cm]{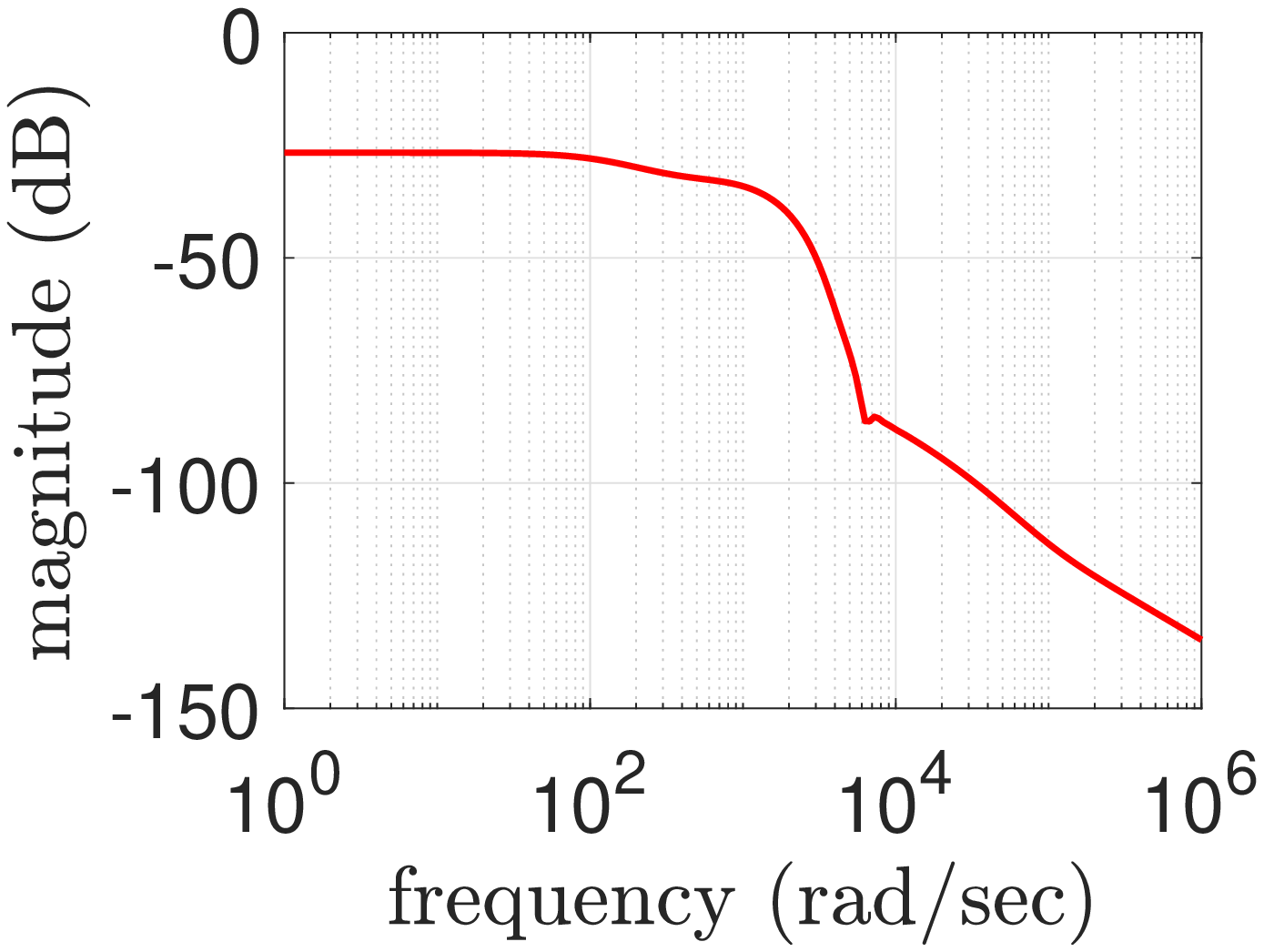}
  \hspace{5mm}
  \includegraphics[width=6.5cm]{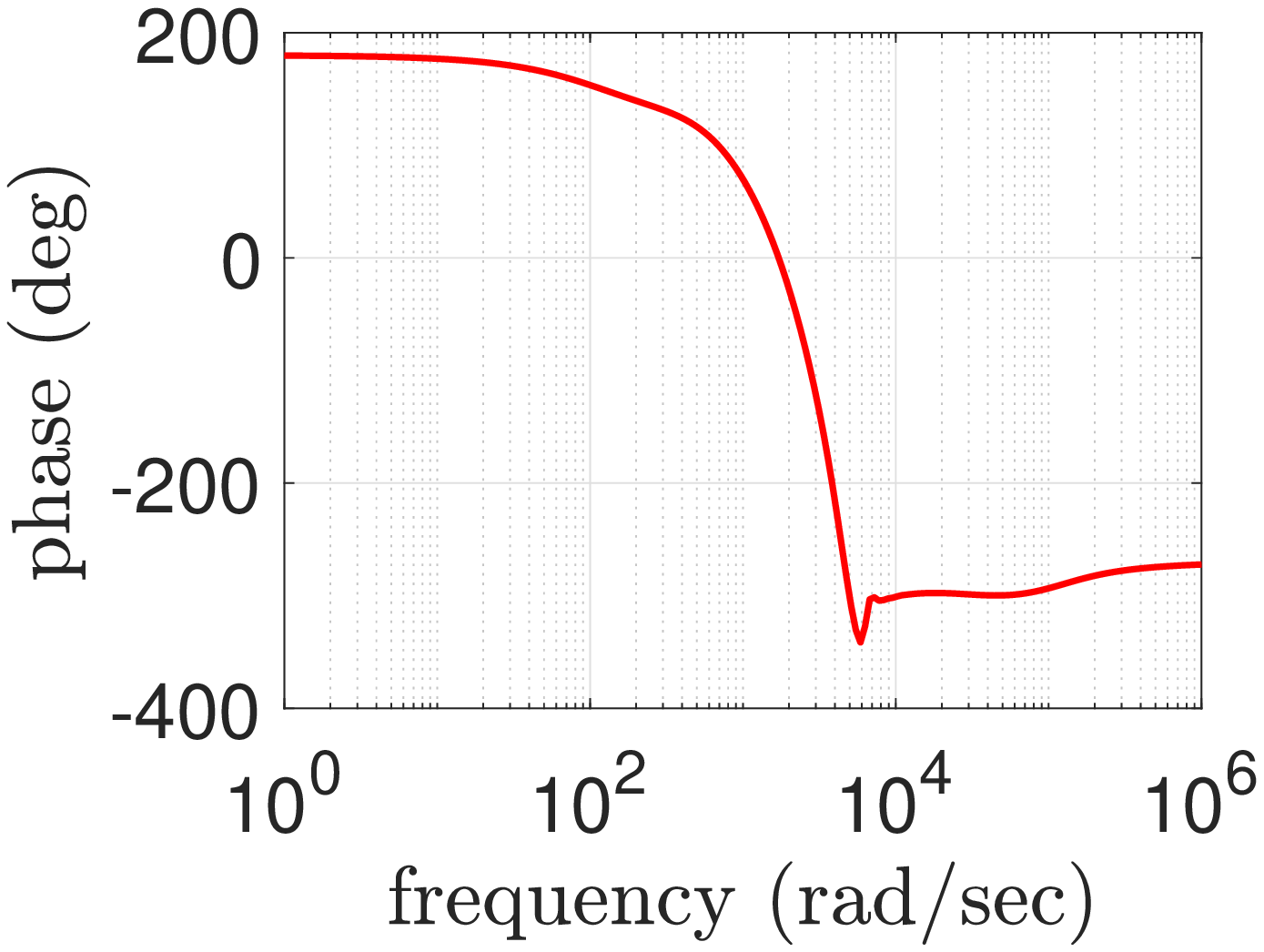}
  \end{center}
  \caption{Bode plot of anemometer benchmark.}
\label{fig:anemometer-bode}
\end{figure}

The anemometer system represents a benchmark from the MOR Wiki~\cite{morwiki}.
A semi-discretisation of a convection-diffusion partial differential
equation yields a linear dynamical system of the form~(\ref{linear-ode})
with SISO. 
Table~\ref{tab:anemometer} depicts its properties.
Since the mass matrix~$E$ is a non-singular diagonal matrix,
we just scale this system to obtain the case of $E=I_n$
used in the following.
Now the matrix~$A$ is stable and non-dissipative.
Figure~\ref{fig:anemometer-bode} illustrates the Bode plot of the system.

\begin{table}
  \caption{Properties of anemometer example. \label{tab:anemometer}}
\begin{center}
  \begin{tabular}{cc} \hline
    dimension~$n$ & 29008 \\
    \# non-zero entries in~$A$ & 201622 \\
    \# non-zero entries in~$E$ & 29008 \\ \hline
  \end{tabular}
\end{center}
\end{table}

\begin{figure}
  \begin{center}
  \includegraphics[width=6.5cm]{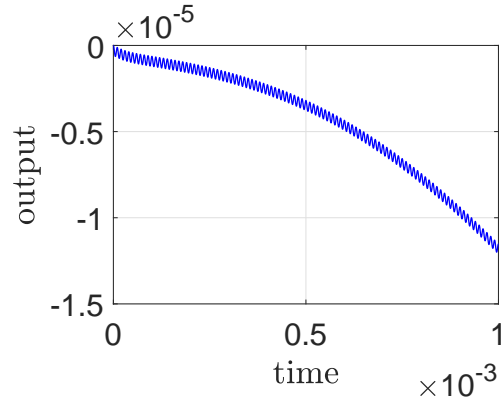}
  \end{center}
  \caption{Output of the anemometer system.}
\label{fig:anemometer-output}
\end{figure}

We perform an MOR by POD as described in~\cite[p.~277]{antoulas}.
The input of the system~(\ref{linear-ode}) is chosen as
the harmonic oscillation
$$ u(t) = \sin \left( \textstyle \frac{2\pi}{T} t \right)
\quad \mbox{with} \quad T = 10^{-5} . $$
We select the initial values $x(0)=0$.
The MATLAB function {\tt ode45} implements a Runge-Kutta method
with time step size selection based on local error control.
This Runge-Kutta method performs 1602 steps and generates $s = 6409$
snapshots (including the internal stages) in the time interval $[0,10^{-3}]$.
The resulting output is depicted in Figure~\ref{fig:anemometer-output}.
We compute a singular value decomposition of the
snapshot matrix $X \in \real^{n \times s}$, where just the largest
singular values are determined.
The $r$ dominant singular vectors build the projection matrix
$V \in \real^{n \times r}$.

\begin{figure}
  \begin{center}
  \includegraphics[width=6.5cm]{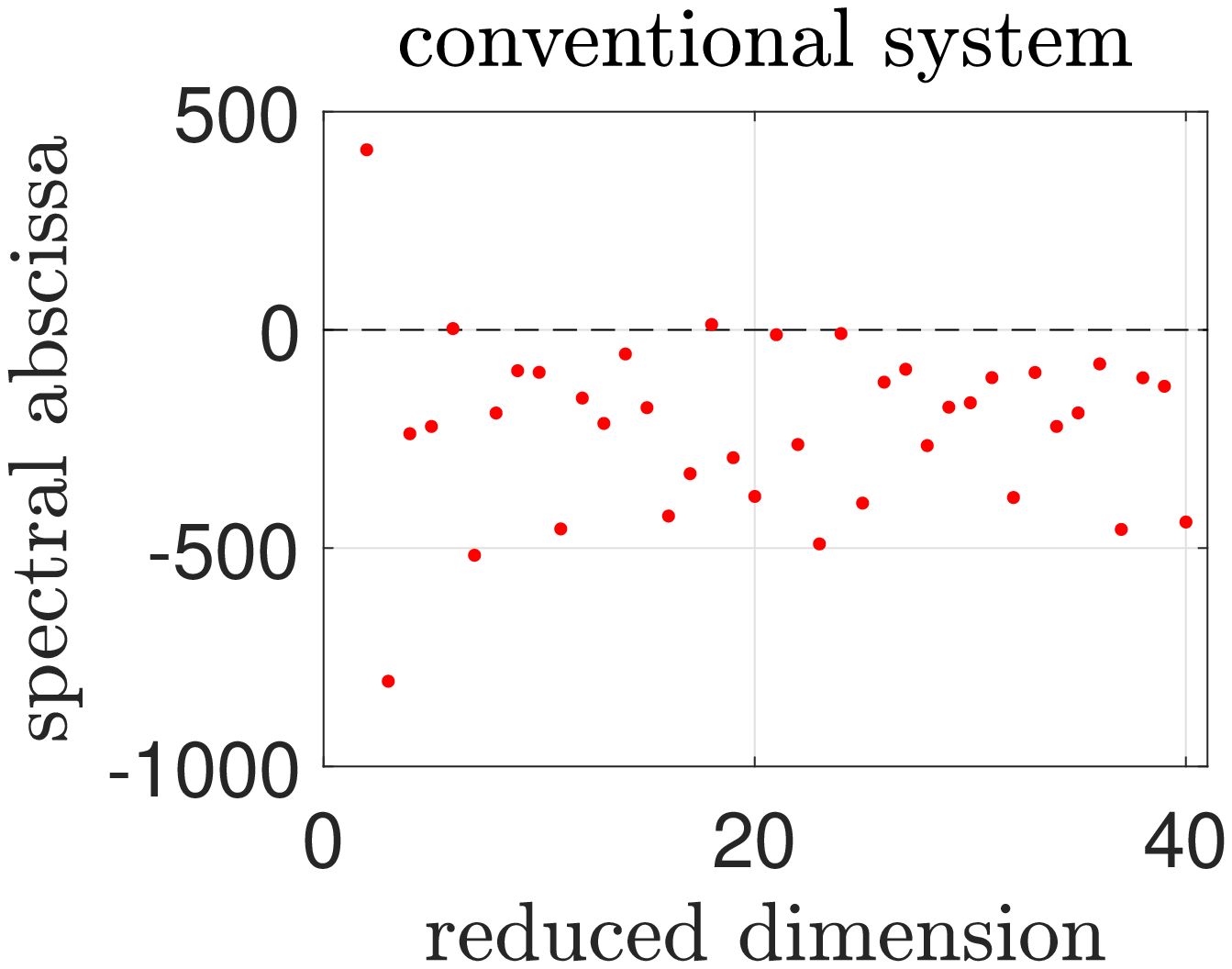}
  \hspace{5mm}
  \includegraphics[width=6.5cm]{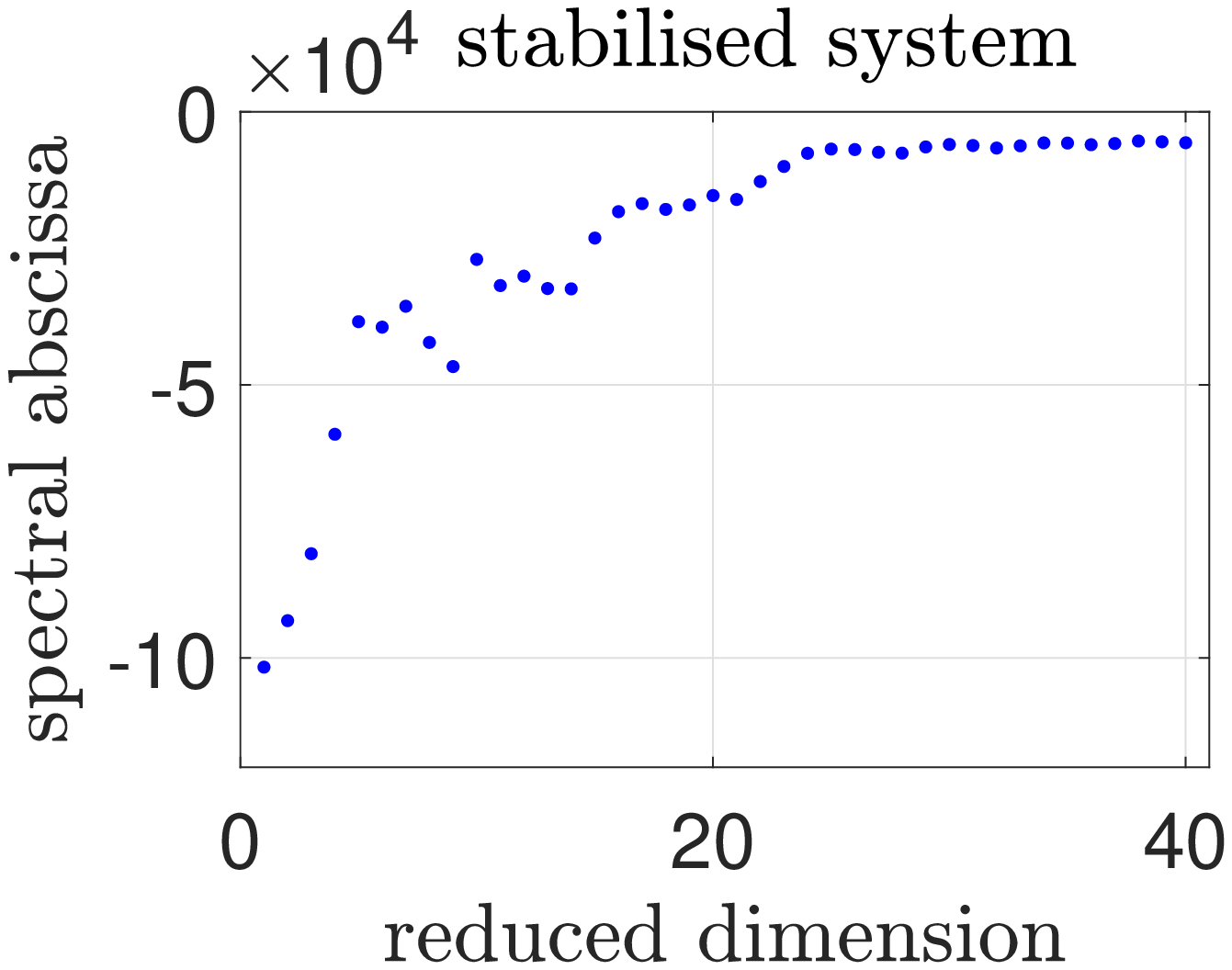}
  \end{center}
  \caption{Spectral abscissa of the matrices in the ROMs
    from conventional system (left)
    and stabilised system (right).}
\label{fig:anemo-spectrum}
\end{figure}

We investigate the cases $r=1,2,\ldots,40$.
The conventional POD method employs the matrix $W=V$
in~(\ref{projected-matrices}).
Figure~\ref{fig:anemo-spectrum} (left) shows the spectral abscissa
from Definition~\ref{def:spectral} for the matrices
$\bar{A}$ ($\bar{E} = I_r$) in the ROMs~(\ref{system-reduced}).
We observe that this MOR generates an unstable system
only in the three cases $r=2,6,18$.

Alternatively, the stability-preserving technique is used
as presented in Section~\ref{sec:lyapunov}.
We compute the largest eigenvalues of the symmetric part $A + A^\top$
($E=I_n$).
It follows that $k = 39$ eigenvalues are non-negative.
Hence the assumption $k \ll n$ is satisfied.
We arrange the matrix~(\ref{delta-F}) with $\delta = 1$ and
solve the Lyapunov equation~(\ref{lyapunov3}) iteratively.
We use the ADI algorithm in the function {\tt lyapchol}
of the sss toolbox, see~\cite{castagnotto}.
Ten iteration steps yield an approximation~(\ref{appr-M}) with
a factor~$Z$ of rank $q=390$.
The projection matrices $V$ and $W$ from~(\ref{matrix-W}) produce
the ROMs~(\ref{system-reduced}), where it holds that $\bar{E} \neq I_r$.
Figure~\ref{fig:anemo-spectrum} (right) illustrates the spectral
abscissa of the matrices $\bar{E}^{-1} \bar{A}$.
All reduced systems are asymptotically stable now.

\begin{figure}
  \begin{center}
  \includegraphics[width=8cm]{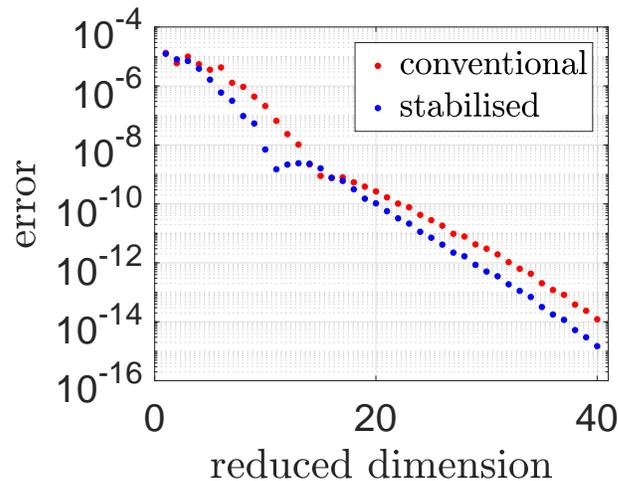}
  \end{center}
  \caption{Maximum error of ROMs for the output in the time domain
    concerning anemometer example.}
\label{fig:anemometer-error}
\end{figure}

We compare the approximation quality of the conventional reduction
and the stabilisation technique.
Numerical solutions of initial value problems for the FOM and the ROMs
are computed by the trapezoidal rule using 1000 time steps of constant size
in $[0,10^{-3}]$.
The maximum absolute errors are depicted for the output
in Figure~\ref{fig:anemometer-error}.
We remind that the magnitude of the output is about $10^{-5}$,
see Figure~\ref{fig:anemometer-output}.
Often the stabilised systems cause a slightly lower error in comparison
to the conventional systems. 
Surprisingly, the instabilities do not affect the
approximation quality of the ROMs in the cases $r=2,6,18$.



\section{Conclusions}
Stability of Galerkin-type projected dynamical systems can be guaranteed
by a transformation, which requires the solution of a Lyapunov equation.
We derived an approximation of the solution, whose unknown part
satisfies an alternative Lyapunov equation including an input matrix
given by a low-rank factorisation.
This approximation ensures a non-singular transformation matrix.
Furthermore, the reduced mass matrices are well-conditioned in general.
Numerical computations confirmed that our stability-preserving
technique is feasible and does not compromise 
the accuracy of the projection-based MOR. 
Moreover, the ADI method produces an efficient iterative solution of the
Lyapunov equation in the case of input matrices with factors
of sufficiently small rank.
However, this restriction on the rank is not always fulfilled.
The rank depends on the number of non-negative eigenvalues
in a symmetric matrix.




\begin{thebibliography}{00}

\bibitem{antoulas}
A.C.~Antoulas,
Approximation of Large-Scale Dynamical Systems, 
SIAM Publications, 2005.

\bibitem{bai-freund}
Z.~Bai, R.~Freund,
A partial Pad{\'e}-via-Lanczos method for reduced order modeling,
Linear Algebra Appl. 332--334 (2001) 139--164.

\bibitem{benner-gugercin-willcox}
P.~Benner, S.~Gugercin, K.~Willcox, 
A survey of projection-based model order reduction methods for
parametric dynamical systems, 
SIAM Review~57 (2015) 483--531.

\bibitem{benner-mehrmann}
P.~Benner, V.~Mehrmann, D.C.~Sorensen,
Dimension Reduction of Large-Scale Systems,
Lecture Notes in Compuational Science and Engineering, Vol.~45,
Springer, 2005.

\bibitem{braun}
M.~Braun,
Differential Equations and Their Applications,
3rd ed., Springer, 1983.

\bibitem{castagnotto}
A.~Castagnotto, M.~Cruz~Varona, L.~Jeschek, B.~Lohmann,
sss \& sssMOR: Analysis and reduction of large-scale dynamic systems in MATLAB,
Auto\-matisierungstechnik~65 (2017) 134--150.

\bibitem{castane-selga}
R.~Casta{\~n}{\'e}~Selga, B.~Lohmann, R.~Eid,
Stability preservation in projection-based model order reduction
of large scale systems,
Eur. J. Control 18 (2012) 122--132.

\bibitem{ding-yao}
J.~Ding, G.~Yao,
The eigenvalue problem of a specially updated matrix,
Appl. Math. Comput. 185 (2007) 415--420.

\bibitem{freund}
R.~Freund, 
Model reduction methods based on Krylov subspaces, 
Acta Numerica 12 (2003) 267--319.

\bibitem{gil}
M.I.~Gil',
Explicit Stability Conditions for Continuous Systems:
A Functional Analytic Approach,
Springer, 2005. 

\bibitem{gugercin-antoulas}
S.~Gugercin, A.C.~Antoulas,
A survey of model reduction by balanced truncation and some new results,
Int. J. Control 77 (2004) 748--766.

\bibitem{hammarling}
S.J.~Hammarling,
Numerical solution of stable non-negative definite Lyapunov equation,
IMA J. Numer. Anal. 2 (1982) 303--323.

\bibitem{ionescu}
T.C.~Ionescu, A.~Astolfi,
On moment matching with preservation of passivity and stability,
in: 49th IEEE Conference on Decision and Control, 2010,
pp. 6189--6194.

\bibitem{kramer}
B.~Kramer, J.R.~Singler,
A POD projection method for large-scale algebraic Riccati equations,
Numer. Algebra Contr. Optim. 6 (2016) 413--435.

\bibitem{li-white}
J.-R.~Li, J.~White,
Low rank solution of Lyapunov equations,
SIAM J. Matrix Anal. \& Appl. 24 (2002) 260--280.

\bibitem{lohmann-eid}
B.~Lohmann, R.~Eid, 
Efficient order reduction of parametric and nonlinear models 
by superposition of locally reduced models, 
in: G.~Roppenecker, B.~Lohmann (eds.), 
Methoden und Anwendungen der Regelungstechnik. Shaker, 2009.

\bibitem{matlab2016}
MATLAB, version 9.1.0.441655 (R2016b), The Mathworks Inc.,
Natick, Massachusetts, 2016.

\bibitem{morwiki}
MOR~Wiki, {\tt https://morwiki.mpi-magdeburg.mpg.de/morwiki} \\
Cited Nov~2, 2017.

\bibitem{mueller}
P.C.~M\"uller,
Modified Lyapunov equations for LTI descriptor systems,
J.~Braz. Soc. Mech. Sci. \& Eng. 28 (2006) 448--452.

\bibitem{penzl}
T.~Penzl, 
LYAPACK: A MATLAB Toolbox for Large Lyapunov and
Riccati Equations, Model Reduction Problems, and
Linear-Quadratic Optimal Control Problems,
Users' Guide (Version 1.0), 1999.

\bibitem{prajna}
S.~Prajna,
POD model reduction with stability guarantee,
in: Proceedings of 42nd IEEE Conference on Decision and Control,
Maui, Hawaii, USA, December 2003, pp. 5254--5258.

\bibitem{pulch18}
R.~Pulch,
Model order reduction and low-dimensional representations for
random linear dynamical systems,
Math. Comput. Simulat. 144 (2018) \mbox{1--20}.

\bibitem{pulch-augustin}
R.~Pulch, F.~Augustin,
Stability preservation in stochastic Galerkin projections of dynamical systems,
arXiv:1708:00958 (2017).

\bibitem{schilders}
W.H.A.~Schilders, M.A.~van~der~Vorst, J.~Rommes (eds.), 
Model Order Reduction: Theory, Research Aspects and Applications, 
Mathematics in Industry, Vol.~13, Springer, 2008.

\bibitem{seydel}
R.~Seydel,
Practical Bifurcation and Stability Analysis,
3rd ed., Springer, 2010.

\bibitem{shmaliy}
Y.~Shmaliy,
Continuous-Time Systems,
Springer, 2007.
  
\bibitem{stoerbulirsch}
J. Stoer, R. Bulirsch,
Introduction to Numerical Analysis,
3rd ed., Springer, New York, 2002.

\bibitem{wolf}
T.~Wolf, H.~Panzer, B.~Lohmann,
Model order reduction by approximate balanced truncation:
a unifying framework,
Automatisierungstechnik 61 (2013) 545--556.

\bibitem{xiu-book}
D.~Xiu,  
Numerical Methods for Stochastic Computations: a Spectral Method Approach, 
Princeton University Press, 2010.

\end{thebibliography}
\end{document}